\font\fifteenrm=cmr10 scaled\magstep2 
\font\fifteeni=cmmi10 scaled\magstep2
\font\fifteensy=cmsy10 scaled\magstep2
\font\fifteenbf=cmbx10 scaled\magstep2
\font\fifteentt=cmtt10 scaled\magstep2
\font\fifteenit=cmti10 scaled\magstep2
\font\fifteensl=cmsl10 scaled\magstep2
\font\fifteenam=msam10 scaled\magstep2
\font\fifteenbm=msbm10 scaled\magstep2
\font\fifteenex=cmex10 scaled\magstep2
\font\fifteensc=cmcsc10 scaled\magstep2 
\font\twelverm=cmr10 at 12pt
\font\twelvei=cmmi10 at 12pt
\font\twelvesy=cmsy10 at 12pt
\font\twelvebf=cmbx10 at 12pt
\font\twelvett=cmtt10 at 12pt
\font\twelveit=cmti10 at 12pt
\font\twelvesl=cmsl10 at 12pt
\font\twelveam=msam10 at 12pt
\font\twelvebm=msbm10 at 12pt
\font\twelveex=cmex10 at 12pt
\font\twelvesc=cmcsc10 at 12pt
\font\elevenrm=cmr10 scaled\magstephalf 
\font\eleveni=cmmi10 scaled\magstephalf
\font\elevensy=cmsy10 scaled\magstephalf
\font\elevenbf=cmbx10 scaled\magstephalf
\font\eleventt=cmtt10 scaled\magstephalf
\font\elevenit=cmti10 scaled\magstephalf
\font\elevensl=cmsl10 scaled\magstephalf
\font\elevenam=msam10 scaled\magstephalf
\font\elevenbm=msbm10 scaled\magstephalf
\font\elevenex=cmex10 scaled\magstephalf
\font\elevensc=cmcsc10 scaled\magstephalf
\font\tenrm=cmr10
\font\teni=cmmi10
\font\tensy=cmsy10
\font\tenbf=cmbx10
\font\tentt=cmtt10
\font\tenit=cmti10
\font\tensl=cmsl10
\font\tenam=msam10
\font\tenbm=msbm10
\font\tenex=cmex10
\font\tensc=cmcsc10
\font\ninerm=cmr9
\font\ninei=cmmi9
\font\ninesy=cmsy9
\font\ninebf=cmbx9
\font\ninett=cmtt9
\font\nineit=cmti9
\font\ninesl=cmsl9
\font\nineam=msam9
\font\ninebm=msbm9
\font\nineex=cmex9
\font\ninesc=cmcsc9
\font\eightrm=cmr8
\font\eighti=cmmi8
\font\eightsy=cmsy8
\font\eightbf=cmbx8
\font\eighttt=cmtt8
\font\eightit=cmti8
\font\eightsl=cmsl8
\font\eightam=msam8
\font\eightbm=msbm8
\font\eightex=cmex8
\font\eightsc=cmcsc8
\font\sevenrm=cmr7
\font\seveni=cmmi7
\font\sevensy=cmsy7
\font\sevenbf=cmbx7

\font\sevenam=msam7
\font\sevenbm=msbm7

\font\sixrm=cmr6
\font\sixi=cmmi6
\font\sixsy=cmsy6

\font\sixam=msam6
\font\sixbm=msbm6

\font\fiverm=cmr5
\font\fivei=cmmi5
\font\fivesy=cmsy5

\font\fiveam=msam5
\font\fivebm=msbm5

\font\fourrm=cmr5 at 4pt
\font\fouri=cmmi5 at 4pt
\font\foursy=cmsy5 at 4pt

\font\fouram=msam5 at 4pt
\font\fourbm=msbm5 at 4pt

\skewchar\twelvei='177 \skewchar\eleveni='177\skewchar\teni='177
\skewchar\ninei='177 \skewchar\eighti='177\skewchar\seveni='177 
\skewchar\sixi='177 \skewchar\fivei='177 \skewchar\fouri='177
\skewchar\twelvesy='60 \skewchar\elevensy='60 \skewchar\tensy='60
\skewchar\ninesy='60 \skewchar\eightsy='60 \skewchar\sevensy='60 
\skewchar\sixsy='60 \skewchar\fivesy='60 \skewchar\foursy='60
\newfam\itfam
\newfam\slfam
\newfam\bffam
\newfam\ttfam
\newfam\scfam
\newfam\amfam
\newfam\bmfam
\def\eightbig#1{{\hbox{$\left#1\vbox to 6.5pt{}\voidright $}}}
\def\eightBig#1{{\hbox{$\left#1\vbox to 7.5pt{}\voidright $}}}
\def\eightbigg#1{{\hbox{$\left#1\vbox to 10pt{}\voidright $}}}
\def\eightBigg#1{{\hbox{$\left#1\vbox to 13pt{}\voidright $}}}
\def\ninebig#1{{\hbox{$\left#1\vbox to 7.5pt{}\voidright $}}}
\def\nineBig#1{{\hbox{$\left#1\vbox to 8.5pt{}\voidright $}}}
\def\ninebigg#1{{\hbox{$\left#1\vbox to 11.5pt{}\voidright $}}}
\def\nineBigg#1{{\hbox{$\left#1\vbox to 14.5pt{}\voidright $}}}
\def\tenbig#1{{\hbox{$\left#1\vbox to 8.5pt{}\voidright $}}}
\def\tenBig#1{{\hbox{$\left#1\vbox to 9.5pt{}\voidright $}}}
\def\tenbigg#1{{\hbox{$\left#1\vbox to 12.5pt{}\voidright $}}}
\def\tenBigg#1{{\hbox{$\left#1\vbox to 16pt{}\voidright $}}}
\def\elevenbig#1{{\hbox{$\left#1\vbox to 9pt{}\voidright $}}}
\def\elevenBig#1{{\hbox{$\left#1\vbox to 10.5pt{}\voidright $}}}
\def\elevenbigg#1{{\hbox{$\left#1\vbox to 14pt{}\voidright $}}}
\def\elevenBigg#1{{\hbox{$\left#1\vbox to 17.5pt{}\voidright $}}}
\def\twelvebig#1{{\hbox{$\left#1\vbox to 10pt{}\voidright $}}}
\def\twelveBig#1{{\hbox{$\left#1\vbox to 11pt{}\voidright $}}}
\def\twelvebigg#1{{\hbox{$\left#1\vbox to 15pt{}\voidright $}}}
\def\twelveBigg#1{{\hbox{$\left#1\vbox to 19pt{}\voidright $}}}
\def\fifteenbig#1{{\hbox{$\left#1\vbox to 12pt{}\voidright $}}}
\def\fifteenBig#1{{\hbox{$\left#1\vbox to 13.5pt{}\voidright $}}}
\def\fifteenbigg#1{{\hbox{$\left#1\vbox to 18pt{}\voidright $}}}
\def\fifteenBigg#1{{\hbox{$\left#1\vbox to 23pt{}\voidright $}}}
\def\voidright{\right.\nulldelimiterspace=0pt \mathsurround=0pt }
\def\fifteenpoint{
  \textfont0=\fifteenrm \scriptfont0=\twelverm \scriptscriptfont0=\tenrm
  \def\rm{\fam0 \fifteenrm}%
  \textfont1=\fifteeni \scriptfont1=\twelvei \scriptscriptfont1=\teni
  \textfont2=\fifteensy \scriptfont2=\twelvesy \scriptscriptfont2=\tensy
  \textfont3=\fifteenex \scriptfont3=\fifteenex \scriptscriptfont3=\fifteenex
  \def\it{\fam\itfam\fifteenit}\textfont\itfam=\fifteenit
  \def\sl{\fam\slfam\fifteensl}\textfont\slfam=\fifteensl
  \def\bf{\fam\bffam\fifteenbf}\textfont\bffam=\fifteenbf 
    \scriptfont\bffam=\twelvebf\scriptscriptfont\bffam=\tenbf
  \def\tt{\fam\ttfam\fifteentt}\textfont\ttfam=\fifteentt
  \def\sc{\fam\scfam\fifteensc}\textfont\scfam=\fifteensc
  \def\am{\fam\amfam\fifteenam}\textfont\amfam=\fifteenam
    \scriptfont\amfam=\twelveam\scriptscriptfont\amfam=\tenam
  \def\bm{\fam\bmfam\fifteenbm}\textfont\bmfam=\fifteenbm
    \scriptfont\bmfam=\twelvebm\scriptscriptfont\bmfam=\tenbm
  \baselineskip=21pt \rm
  \let\big=\fifteenbig\let\Big=\fifteenBig\let\bigg=\fifteenbigg
  \let\Bigg=\fifteenBigg}
\def\twelvepoint{
  \textfont0=\twelverm \scriptfont0=\ninerm \scriptscriptfont0=\sevenrm
  \def\rm{\fam0 \twelverm}%
  \textfont1=\twelvei \scriptfont1=\ninei \scriptscriptfont1=\seveni
  \textfont2=\twelvesy \scriptfont2=\ninesy \scriptscriptfont2=\sevensy
  \textfont3=\twelveex \scriptfont3=\twelveex \scriptscriptfont3=\twelveex
  \def\it{\fam\itfam\twelveit}\textfont\itfam=\twelveit
  \def\sl{\fam\slfam\twelvesl}\textfont\slfam=\twelvesl
  \def\bf{\fam\bffam\twelvebf}\textfont\bffam=\twelvebf 
    \scriptfont\bffam=\ninebf\scriptscriptfont\bffam=\sevenbf
  \def\tt{\fam\ttfam\twelvett}\textfont\ttfam=\twelvett
  \def\sc{\fam\scfam\twelvesc}\textfont\scfam=\twelvesc
  \def\am{\fam\amfam\twelveam}\textfont\amfam=\twelveam
    \scriptfont\amfam=\nineam\scriptscriptfont\amfam=\sevenam
  \def\bm{\fam\bmfam\twelvebm}\textfont\bmfam=\twelvebm
    \scriptfont\bmfam=\ninebm\scriptscriptfont\bmfam=\sevenbm
  \baselineskip=17.8pt \rm 
  \def\looselineskip{\baselineskip=18.5pt plus 1.8pt}%
  \def\tightlineskip{\baselineskip=16.5pt}%
  \def\verytightlineskip{\baselineskip=15pt}%
  \let\big=\twelvebig\let\Big=\twelveBig\let\bigg=\twelvebigg
  \let\Bigg=\twelveBigg  }
\def\elevenpoint{
  \textfont0=\elevenrm \scriptfont0=\ninerm \scriptscriptfont0=\sixrm
  \def\rm{\fam0 \elevenrm}%
  \textfont1=\eleveni \scriptfont1=\ninei \scriptscriptfont1=\sixi
  \textfont2=\elevensy \scriptfont2=\ninesy \scriptfont2=\sixsy 
  \textfont3=\elevenex \scriptfont3=\elevenex \scriptfont3=\elevenex
  \def\it{\fam\itfam\elevenit}\textfont\itfam=\elevenit
  \def\sl{\fam\slfam\elevensl}\textfont\slfam=\elevensl
  \def\bf{\fam\bffam\elevenbf}\textfont\bffam=\elevenbf
  \def\tt{\fam\ttfam\eleventt}\textfont\ttfam=\eleventt
  \def\sc{\fam\scfam\elevensc}\textfont\scfam=\elevensc
  \def\am{\fam\amfam\elevenam}\textfont\amfam=\elevenam
    \scriptfont\amfam=\nineam\scriptscriptfont\amfam=\sixam
  \def\bm{\fam\bmfam\elevenbm}\textfont\bmfam=\elevenbm
    \scriptfont\bmfam=\ninebm\scriptscriptfont\bmfam=\sixbm
  \baselineskip=15.1pt \rm
  \def\looselineskip{\baselineskip=16pt plus 1.5pt}%
  \def\tightlineskip{\baselineskip=14pt}%
  \def\verytightlineskip{\baselineskip=13pt}%
  \let\big=\elevenbig\let\Big=\elevenBig\let\bigg=\elevenbigg
  \let\Bigg=\elevenBigg  }
\def\tenpoint{
  \textfont0=\tenrm \scriptfont0=\eightrm \scriptscriptfont0=\fiverm
  \def\rm{\fam0 \tenrm}%
  \textfont1=\teni \scriptfont1=\eighti \scriptscriptfont1=\fivei
  \textfont2=\tensy \scriptfont2=\eightsy \scriptfont2=\fivesy 
  \textfont3=\tenex \scriptfont3=\tenex \scriptfont3=\tenex
  \def\it{\fam\itfam\tenit}\textfont\itfam=\tenit
  \def\sl{\fam\slfam\tensl}\textfont\slfam=\tensl
  \def\bf{\fam\bffam\tenbf}\textfont\bffam=\tenbf
  \def\tt{\fam\ttfam\tentt}\textfont\ttfam=\tentt
  \def\sc{\fam\scfam\tensc}\textfont\scfam=\tensc
  \def\am{\fam\amfam\tenam}\textfont\amfam=\tenam
    \scriptfont\amfam=\eightam \scriptscriptfont\amfam=\fiveam
  \def\bm{\fam\bmfam\tenbm}\textfont\bmfam=\tenbm
    \scriptfont\bmfam=\eightbm \scriptscriptfont\bmfam=\fivebm
  \baselineskip=14pt\rm
  \def\looselineskip{\baselineskip=14.8pt plus1.5pt}
  \def\tightlineskip{\baselineskip=12.6pt}%
  \def\verytightlineskip{\baselineskip=13pt}%
  \let\big=\tenbig\let\Big=\tenBig\let\bigg=\tenbigg\let\Bigg=\tenBigg  }
\def\ninepoint{
  \textfont0=\ninerm \scriptfont0=\sevenrm \scriptscriptfont0=\fourrm
  \def\rm{\fam0 \ninerm}%
  \textfont1=\ninei \scriptfont1=\seveni \scriptscriptfont1=\fouri
  \textfont2=\ninesy \scriptfont2=\sevensy \scriptfont2=\foursy 
  \textfont3=\nineex \scriptfont3=\nineex \scriptfont3=\nineex
  \def\it{\fam\itfam\nineit}\textfont\itfam=\nineit
  \def\sl{\fam\slfam\ninesl}\textfont\slfam=\ninesl
  \def\bf{\fam\bffam\ninebf}\textfont\bffam=\ninebf
  \def\tt{\fam\ttfam\ninett}\textfont\ttfam=\ninett
  \def\sc{\fam\scfam\ninesc}\textfont\scfam=\ninesc
  \def\am{\fam\amfam\nineam}\textfont\amfam=\nineam
    \scriptfont\amfam=\nineam\scriptscriptfont\amfam=\fouram
  \def\bm{\fam\bmfam\ninebm}\textfont\bmfam=\ninebm
    \scriptfont\bmfam=\ninebm\scriptscriptfont\bmfam=\fourbm
  \baselineskip=12.6pt\rm
  \def\tightlineskip{\baselineskip=11.5pt}
  \let\big=\ninebig\let\Big=\nineBig\let\bigg=\ninebigg
  \let\Bigg=\nineBigg  }
\def\eightpoint{
  \textfont0=\eightrm \scriptfont0=\fiverm \scriptscriptfont0=\fourrm
  \def\rm{\fam0 \eightrm}%
  \textfont1=\eighti \scriptfont1=\fivei \scriptscriptfont1=\fouri
  \textfont2=\eightsy \scriptfont2=\fivesy \scriptfont2=\foursy 
  \textfont3=\eightex \scriptfont3=\eightex \scriptfont3=\eightex
  \def\it{\fam\itfam\eightit}\textfont\itfam=\eightit
  \def\sl{\fam\slfam\eightsl}\textfont\slfam=\eightsl
  \def\bf{\fam\bffam\eightbf}\textfont\bffam=\eightbf
  \def\tt{\fam\ttfam\eighttt}\textfont\ttfam=\eighttt
  \def\sc{\fam\scfam\eightsc}\textfont\scfam=\eightsc
  \def\am{\fam\amfam\eightam}\textfont\amfam=\eightam
    \scriptfont\amfam=\eightam\scriptscriptfont\amfam=\fouram
  \def\bm{\fam\bmfam\eightbm}\textfont\bmfam=\eightbm
    \scriptfont\bmfam=\eightbm\scriptscriptfont\bmfam=\fourbm
  \baselineskip=11.2pt \rm
  \def\tightlineskip{\baselineskip=10.4pt}
  \let\big=\eightbig\let\Big=\eightBig\let\bigg=\eightbigg
  \let\Bigg=\eightBigg  }

\twelvepoint
%
\twelvepoint
\nopagenumbers
\hsize=6in\vsize=8.8in

\parskip=1pt plus 1pt

\newif\ifSpecialhead\Specialheadfalse
\newbox\specialheadbox

\def\specialhead #1\par{\Specialheadtrue\setbox\specialheadbox=\hbox{#1}}
\headline={{\ifSpecialhead\box\specialheadbox\global\Specialheadfalse\else
     \ifnum\pageno<0{\hfill\quad{\twelvebf\folio}}%
     \else\ifnum\pageno<2\hfill
     \else\hfill\twelvepoint\sc\firstmark\quad{\twelvebf\folio}\fi\fi\fi}}

\def\title#1\par{\medskip{\def\cr{\hfil\par\hfil}\hfil\fifteenbf#1\hfil\par}\medskip}
\def\subtitle#1\par{\centerline{\fifteenrm #1}\medskip}
\def\author#1\par{\medskip{\def\cr{\hfil\par\hfil\twelvesc}\fifteensc\hfil#1\hfil\par}}
\def\authors#1\par{\medskip\fifteensc\center#1\par}
\def\center#1\par{{\def\cr{\hfil\par\hfil}\hfil#1\hfil\par}}
\def\abstract.#1\par{\message{Abstract.}%
                    \medskip{\narrower\narrower\tenpoint\tightlineskip
                        \noindent{\bf Abstract.}#1\par}\medskip\noindent}
\def\tinyabstract.#1\par{\message{Abstract.}%
                    \medskip{\narrower\narrower\eightpoint\tightlineskip
                        \noindent{\bf Abstract.}#1\par}\medskip\noindent}
\def\bigabstract.#1\par{\message{Abstract.}%
                         \medskip{\narrower\narrower\tightlineskip
                         \noindent{\bf Abstract. }#1\par}\medskip\noindent}
\def\acknowledgement#1\par{\footnote{}{#1}}
\def\sectionskip{\Goodbreak\vskip 25pt plus 15pt minus 5pt}
\def\secnumber{\ifquiet
               \else\ifNoSections
                    \else\sectionsymbol\the\secno\quad\fi\fi}
\def\section#1\par{ \NoSectionsfalse\par\sectionskip\proofdepth=0\claimno=0
 \ifquiet\else\advance\secno by1\fi\toks0={#1}
 \immediate\write16{\ifquiet\else Section \the\secno\space\fi
                    \the\toks0}%
 \mark{\secnumber #1}%
 {\fifteenpoint\bf\noindent\secnumber #1}\nobreak\bigskip\quietoff
 \nobreak\noindent}

\def\QUIET{\QUIETtrue\quiettrue}

\def\quietoff{\ifQUIET\else\quietfalse\fi}
\newif\ifquiet
\newif\ifQUIET
\newif\ifNoSections
\newcount\claimtype
\newcount\secno
\newcount\claimno
\newcount\subclaimno
\newcount\subsubclaimno
\newcount\subsubsubclaimno
\newcount\proofdepth
\def\subclaimnumber{\ifquiet\else\ifcase\subclaimno\or A\or B\or C\or D\or E\or
     F\or G\or H\or I\or J\or K\or L\or M\or N\or O\or P\fi\fi}
\def\subsubclaimnumber{\ifquiet\else\ifcase\subsubclaimno\or i\or ii\or iii\or
   iv\or v\or vi\or vii\or viii\or ix\or x\or xi\or xii\or xiii\or xiv\fi\fi}
\def\subsubsubclaimnumber{\ifquiet\else\ifcase\subsubsubclaimno\or a\or b\or
   c\or d\or e\or f\or g\or \or h\or i\or j\or k\or l\or m\or n\or o\fi\fi}
\def\claimtag{\ifquiet\else
  \ifNoSections
    \ifcase\proofdepth\the\claimno%
    \or\the\claimno.\subclaimnumber
    \or\the\claimno.\subclaimnumber.\subsubclaimnumber
    \or\the\claimno.\subclaimnumber.\subsubclaimnumber
                                                .\subsubsubclaimnumber\fi
  \else
    \ifcase\proofdepth\the\secno.\the\claimno
    \or\the\secno.\the\claimno.\subclaimnumber
    \or\the\secno.\the\claimno.\subclaimnumber.\subsubclaimnumber
    \or\the\secno.\the\claimno.\subclaimnumber.\subsubclaimnumber
                                                .\subsubsubclaimnumber\fi\fi\fi}
\secno=0\claimno=0\proofdepth=0\subclaimno=0\subsubclaimno=0\subsubsubclaimno=0
\NoSectionstrue
\newbox\qedbox
\def\claimname{\ifcase\claimtype Theorem\or Lemma\or Claim\or Corollary\or
               Question\or Definition\or Remark\or Conjecture\fi}
\def\preclaimskip{\removelastskip
    \ifcase\claimtype\goodbreak\vskip 8pt plus 4pt minus 2pt
                  \or\goodbreak\vskip 6pt plus 4pt minus 1pt
                  \or\goodbreak\vskip 5pt plus 4pt minus 1pt
                  \or\goodbreak\vskip 8pt plus 4pt minus 2pt
                  \or\vskip 7pt plus 4pt minus 2pt
                  \or\vskip 7pt plus 4pt minus 2pt
                  \or\vskip 7pt plus 4pt minus 2pt
                  \or\goodbreak\vskip 8pt plus 4pt minus 2pt\fi}
\def\postclaimskip{\ifcase\claimtype         \vskip 4pt plus 2pt minus 2pt
                                          \or\vskip 3pt plus 2pt minus 2pt
                                          \or\vskip 2pt plus 2pt minus 1pt
                                          \or\vskip 4pt plus 2pt minus 2pt
                                          \or\vskip 1pt plus 2pt
                                          \or\vskip 4pt plus 4pt
                                          \or\vskip 3pt plus 2pt
                                          \or\vskip 4pt plus 2pt minus 2pt\fi}
\def\claimfont{\ifcase\claimtype
                  \sl\or\sl\or\sl\or\sl\or\sl\or\rm\or\rm\or\sl\fi}
\def\advancetag{\ifcase\proofdepth\advance\claimno by1
                               \or\advance\subclaimno by1
                               \or\advance\subsubclaimno by1
                               \or\advance\subsubsubclaimno by1\fi}
\def\sayclaim#1.#2 #3\par{\ifquiet\else\advancetag\fi
    \preclaimskip\setbox1=\hbox{#1}\setbox2=\hbox{#2}%
    \toks0={#1 }
    \immediate\write16{\ifdim\wd1>0pt\the\toks0
                       \else\claimname\space\fi \claimtag.}%
    \vbox{\noindent
    {\bf\ifdim\wd1=0pt \claimname\else #1\fi\ifquiet.\else\ \claimtag{\ifNoSections.\fi}\fi}%
    \enspace{\ifdim\wd2>0pt\sc #2\enspace\fi}%
    {\claimfont #3\par}}\postclaimskip\quietoff}
\def\theorem{\claimtype=0\sayclaim}
\def\lemma{\claimtype=1\sayclaim}
\def\claim{\claimtype=2\sayclaim}
\def\corollary{\claimtype=3\sayclaim}
\def\question{\claimtype=4\sayclaim}

\def\point#1. #2\par{\item{\rm #1.}#2\par}
\def\points#1\cr\par{\medskip\vbox{\let\cr=\point\point#1\par}\par}
\def\df{\it}
\def\prooffont{}
\def\proofsize{}
\def\proofindent{}
\def\proofskip{\badbreak\ifcase\claimtype    \vskip 3pt plus 2pt minus 2pt
                                          \or\vskip 2pt plus 2pt minus 2pt
                                          \or\vskip 1pt plus 2pt minus 1pt
                                          \or\vskip 3pt plus 2pt minus 2pt
                                          \or\vskip 1pt plus 2pt
                                          \or\vskip 2pt plus 4pt
                                          \or\vskip 1pt plus 2pt
                                          \or\vskip 3pt plus 2pt minus 2pt\fi}

\def\Goodbreak{\vskip0pt plus.5in\penalty-1000\vskip0pt plus-.5in}
\def\goodbreak{\penalty-500}
\def\badbreak{\penalty500}
\def\Badbreak{\penalty1000}
\def\proof{\message{proof}\removelastskip\Badbreak\proofskip\begingroup
  \advance\proofdepth by1
  \setbox\qedbox=\hbox{\halmos\raise2pt\hbox{\fiverm\claimname}}%
  \prooffont\proofsize\proofindent\noindent{\bf Proof: }}
\def\proofof#1:{\message{proof}\removelastskip\Badbreak\proofskip\begingroup
  \advance\proofdepth by1
  \setbox\qedbox=\hbox{\halmos\raise2pt\hbox{\fiverm#1}}%
  \prooffont\proofsize\proofindent\noindent{\bf Proof of #1: }}
\def\cite[#1]{[{\tenrm{#1}}]\message{[#1]}}
\edef\ref#1{\expandafter\global\expandafter\edef#1{\noexpand\claimtag}}
\newwrite\notes
\openout\notes=\jobname.notes
\long\def\unexpandedwrite#1#2{\def\finwrite{\write#1}%
   {\aftergroup\finwrite\aftergroup{\sanitize#2\endsanity}}}
\def\sanitize{\futurelet\next\sanswitch}
\let\stoken=\space
\def\sanswitch{\ifx\next\endsanity
  \else\ifcat\noexpand\next\stoken\aftergroup\space\let\next=\eat
   \else\ifcat\noexpand\next\bgroup\aftergroup{\let\next=\eat
    \else\ifcat\noexpand\next\egroup\aftergroup}\let\next=\eat
     \else\let\next=\copytoken\fi\fi\fi\fi \next}
\def\eat{\afterassignment\sanitize \let\next= }
\long\def\copytoken#1{\ifcat\noexpand#1\relax\aftergroup\noexpand
  \else\ifcat\noexpand#1\noexpand~\aftergroup\noexpand\fi\fi
  \aftergroup#1\sanitize}
\def\endsanity\endsanity{}

\def\note#1#2{\hbox to2in{\strut#1\quad\dotfill\quad#2}}
\def\boxit#1{\setbox4=\hbox{\kern1pt#1\kern1pt}
  \hbox{\vrule\vbox{\hrule\kern1pt\box4\kern1pt\hrule}\vrule}}
\def\halmos{\hbox{\am\char'3}}
\def\qed#1\par{\message{.                                }\setbox1=\hbox{#1}%
  \ifdim\wd1>0pt\setbox\qedbox=\hbox{\halmos\raise2pt\hbox{\fiverm #1}}\fi
  \kern5pt\lower 2pt\hbox{\box\qedbox}\proofskip\goodbreak\endgroup}

\def\sectionsymbol{\S}
\def\k{\kappa}
\def\g{\gamma}
\def\a{\alpha}
\def\b{\beta}
\def\d{\delta}
\def\s{\sigma}
\def\t{\tau}
\def\l{\lambda}

\def\I1{\mathop{\hbox{\sc i}_1}}

\def\P{{\mathchoice{\hbox{\bm P}}{\hbox{\bm P}}
         {\hbox{\tenbm P}}{\hbox{\sevenbm P}}}}
\def\Q{{\mathchoice{\hbox{\bm Q}}{\hbox{\bm Q}}
         {\hbox{\tenbm Q}}{\hbox{\sevenbm Q}}}}
\def\R{{\mathchoice{\hbox{\bm R}}{\hbox{\bm R}}
         {\hbox{\tenbm R}}{\hbox{\sevenbm R}}}}

\def\card#1{\left|#1\right|}

\def\dom{\mathop{\rm dom}\nolimits}

\def\elesub{\prec}

\def\unifto{\buildrel\lower 7pt\hbox{$\to$}\over\to}

\def\cof{\mathop{\rm cof}\nolimits}
\def\cp{\mathop{\rm cp}\nolimits}

\def\plus{^{\scriptscriptstyle +}}
\def\plusplus{^{\scriptscriptstyle ++}}

\def\in{\mathrel{\mathchoice{\raise
1pt\hbox{$\scriptstyle\cal\char'62$}}
         {\raise 1pt\hbox{$\scriptstyle\cal\char'62$}}
         {\raise .5pt\hbox{$\scriptscriptstyle\cal\char'62$}}
         {\hbox{$\scriptscriptstyle\cal\char'62$}}}\penalty700{}}
\def\ni{\mathrel{\mathchoice{\raise 1pt\hbox{$\scriptstyle\cal\char'63$}}
                   {\raise 1pt\hbox{$\scriptstyle\cal\char'63$}}
                   {\raise .5pt\hbox{$\scriptscriptstyle\cal\char'63$}}
                   {\hbox{$\scriptscriptstyle\cal\char'63$}}}\penalty700}
\def\of{\mathrel{\mathchoice{\raise 1pt\hbox{$\scriptstyle\subseteq$}}
                   {\raise 1pt\hbox{$\scriptstyle\subseteq$}}
                   {\raise .5pt\hbox{$\scriptscriptstyle\subseteq$}}
                   {\hbox{$\scriptscriptstyle\subseteq$}}}}
\def\fo{\mathrel{\mathchoice{\raise 1pt\hbox{$\scriptstyle\supseteq$}}
                   {\raise 1pt\hbox{$\scriptstyle\supseteq$}}
                   {\raise .5pt\hbox{$\scriptscriptstyle\supseteq$}}
                   {\hbox{$\scriptscriptstyle\supseteq$}}}}
\def\notin{\mathrel{\mathchoice
  {\raise 1pt\hbox{\rlap{$\scriptstyle\;|$}$\scriptstyle\cal\char'62$}}
  {\raise 1pt\hbox{\rlap{$\scriptstyle\kern2pt
          |$}$\scriptstyle\cal\char'62$}}
  {\raise .5pt\hbox{\rlap{$\scriptscriptstyle\, |$}$\scriptscriptstyle
      \cal\char'62$}}
  {\hbox{\rlap{$\scriptscriptstyle\, |$}$\scriptscriptstyle
     \cal\char'62$}}}%
  \penalty700}

\def\ofnoteq{\mathbin{\hbox{\bm\char'050}}}

\def\and{\mathrel{\kern1pt\&\kern1pt}}

\def\implies{\rightarrow}

\def\union{\cup}

\def\intersect{\cap}

\def\setminus{\mathbin{\hbox{\bm\char'162}}}

\def\add{\mathop{\rm Add}\nolimits}

\def\cross{\times}

\def\muchgt{\gg}
\def\lt{\mathrel{\mathchoice{\scriptstyle<}{\scriptstyle<}
   {\scriptscriptstyle<}{\scriptscriptstyle<}}}
\def\lte{\mathrel{\scriptstyle\leq}}

\def\[#1]{\left[\vphantom{\bigm|}#1\right]}
\def\<#1>{\langle\,#1\,\rangle}

\def\restrict{\mathbin{\mathchoice{\hbox{\am\char'26}}{\hbox{\am\char'26}}{\hbox{\eightam\char'26}}{\hbox{\sixam\char'26}}}}
\def\force{\mathbin{\hbox{\am\char'15}}}

\def\emptyset{\mathord{\hbox{\bm\char'77}}}

\def\boolval#1{\mathopen{\lbrack\!\lbrack}\,#1\,\mathclose{\rbrack\!\rbrack}}
\def\concat{\mathbin{{\,\hat{ }\,}}}

\def\decides{\mathrel{{||}}}

\def\st{\mid}
\def\seq<#1>{{\def\st{\mid\penalty650}\left<\,#1\,\right>}}

\def\set#1{\{\,{#1}\,\}}

\def\forces{\force}
\def\lttheta{{\raise 1pt\hbox{$\scriptstyle<$}\theta}}

\def\I1{\mathop{\hbox{\sc i}_1}}
\def\ltk{{{\scriptstyle<}\k}}

\def\ltd{{{\scriptstyle<}\d}}

\def\lteb{{{\scriptstyle\leq}\b}}
\def\lted{{{\scriptstyle\leq}\d}}

\def\Rdot{\dot\R}
\def\Qdot{\dot\Q}
\def\PQ{{\P*\Qdot}}
\def\Adot{\dot A}
\def\Qdot{\dot\Q}
\def\qdot{\dot q}

\def\Pforces{\force_{\P}}
\def\PQforces{\force_{\PQ}\;}
\def\Qforces{\force_{\Q}}
\def\Qlte{\lte_{\Q}}
\def\Qterm{\Q_{\rm\scriptscriptstyle term}}
\def\Gterm{G_{\rm\scriptscriptstyle term}}
\def\Ddot{\dot D}
\def\Dterm{D_{\rm\scriptscriptstyle term}}

\def\hdot{\dot h}
\def\rdot{\dot r}
\def\sdot{\dot s}
\def\pbar{\bar p}

\def\one{1\hskip-3pt {\rm l}}

\def\Qterm{\Q_{\fiverm term}}

\def\qdot{\dot q}
\def\rdot{\dot r}
\def\sdot{\dot s}

\def\forces{\force}
\def\I1{\mathop{\hbox{\sc i}_1}}
\def\ltk{{{\scriptstyle<}\k}}
\def\ltd{{{\scriptstyle<}\d}}

\def\lteb{{{\scriptstyle\leq}\b}}
\def\lted{{{\scriptstyle\leq}\d}}

\def\k{\kappa}
\def\g{\gamma}
\def\a{\alpha}
\def\b{\beta}
\def\d{\delta}
\def\s{\sigma}
\def\t{\tau}
\def\l{\lambda}
\def\Rdot{\dot\R}
\def\Qdot{\dot\Q}
\def\Rdd{\ddot\R}
\def\PQ{{\P*\Qdot}}
\def\pq{{\<p,\dot q>}}
\def\pr{{\<p,\dot r>}}
\def\ps{{\<p,\dot s>}}
\def\Adot{\dot A}
\def\Qdot{\dot\Q}
\def\qdot{\dot q}

\def\Pforces{\force_{\P}}
\def\PQforces{\force_{\PQ}\;}
\def\Qforces{\force_{\Q}}
\def\Qlte{\lte_{\Q}}
\def\Qterm{\Q_{\rm\scriptscriptstyle term}}
\def\Gterm{G_{\rm\scriptscriptstyle term}}
\def\Ddot{\dot D}
\def\Dterm{D_{\rm\scriptscriptstyle term}}

\def\termlte{\lte_{\rm\scriptscriptstyle term}}
\def\one{1}
\def\hdot{\dot h}
\def\qdot{\dot q}
\def\rdot{\dot r}
\def\sdot{\dot s}
\def\pgqg{\<p_{\g},\qdot_{\g}>}

\def\pbar{\bar p}
\def\gG{g*G}
\def\VgG{V[g][G]}
\QUIET

\title Small Forcing Makes Any Cardinal Superdestructible

\author Joel David Hamkins

\abstract. Small forcing always ruins the indestructibility of
an indestructible supercompact cardinal. In fact, after small forcing, any
cardinal $\k$ becomes {\it superdestructible}---any
further $\ltk$-closed forcing which adds a subset to $\k$ will
destroy the measurability, even the weak compactness, of $\k$.
Nevertheless, after small
forcing indestructible cardinals remain resurrectible, but never
strongly resurrectible.

Arthur Apter, motivated by issues arising in his recent paper
\cite[AS] with Saharon Shelah, asked me the following
question: ``Does small forcing preserve the
indestructibility of a supercompact cardinal after the Laver
preparation?'' While it is tempting to believe that all large
cardinal properties are preserved by small forcing, the fact
is that the answer to his question is no. Even adding a
Cohen real ruins the indestructibility of any
cardinal. What's more, it is ruined in a very strong way. In this
paper I will prove that small forcing makes any cardinal superdestructible.

Before stating my theorem, let me make some definitions.
In one of my favorite arguments, Laver \cite[L]
proved that with the proper preparation, now called the Laver
preparation, a supercompact cardinal $\k$ can be
made {\df indestructible} in the sense that any $\ltk$-directed closed
forcing preserves the supercompactness of $\k$.
We say that $\k$ is {\df destructible}, therefore, when some
$\ltk$-directed
closed poset destroys the supercompactness of $\k$. Going beyond this,
define that $\k$ is {\df superdestructible} when every
$\ltk$-closed forcing which adds a subset to $\k$ destroys the
measurability of $\k$, and that $\k$ is {\df superdestructible at $\l$}
when any $\ltk$-closed forcing which adds a subset to $\l$ destroys
the $\l$-supercompactness of $\k$. Define $\k$ to be
{\df resurrectible} iff whenever a $\ltk$-directed closed forcing
$\Q$ happens to destroy the supercompactness of $\k$, it can
nevertheless be restored with further $\ltk$-distributive forcing
$\Rdot$; and $\k$ is {\df strongly resurrectible} when $\Rdot$ can be made
actually $\ltk$-closed (this resembles the notion for huge
cardinals in \cite[B]).
Finally, a poset $\P$ is {\df small} relative to
$\k$ when $\card{\P}<\k$. Throughout I consider only nontrivial
posets---forcing with them must add some new set.
Now I am ready to state my main theorem.

\theorem Main Theorem. Small forcing makes any cardinal
superdestructible. Indeed, after small forcing, any
$\ltk$-closed forcing which adds a subset to $\k$ will destroy
the weak compactness of $\k$. What's more, after small forcing, $\k$
becomes superdestructible at $\k\plus$, $\k\plusplus$, etc.
Nevertheless, after small forcing an indestructible cardinal remains
resurrectible, but never strongly resurrectible.

I will actually prove a better theorem: after forcing of size $\b<\k$,
any $\lteb$-closed forcing which adds a subset but no bounded subset
to $\k$ will destroy the measurability and weak compactness of $\k$.
After adding a Cohen real, for example, any countably closed poset
which adds a subset but no bounded subset to $\k$ will destroy the
measurability of $\k$.

This theorem is related to my Fragile Measurability theorems in \cite[H].
There, I show how to force from a model in which $\k$ is
strong, supercompact, or $\I1$, while preserving this large
cardinal property, to a model in which the measurability of
$\k$ is {\df fragile} in the sense that it is destroyed by
any forcing which preserves $\k^{\ltk}$, $\k\plus$, but not $P(\k)$.
To get superdestructibility from fragility we drop the requirement
that $\k\plus$ is preserved, but require the poset to be a little closed.
The two properties are similar in that if $\k$ is
fragile or superdestructible, the measurability of $\k$ is easily
destroyed by forcing. In my fragile measurability
models \cite[H], $\k$ is both fragile and superdestructible.

What is perhaps the first theorem in this line is due to W. Hugh Woodin
\cite[W],
who forced to a model of a supercompact cardinal $\k$ whose measurability
and weak compactness is destroyed by the poset $\add(\k,1)=\k^{\ltk}$.
Woodin used a
reverse Easton $\k$-iteration, adding a system of coherent clubs. Later,
he simplified his argument to add just a subset of $\d$ at certain
stages $\d$. My theorems here show that the entire $\k$-iteration may be
replaced by any small forcing, such as adding a Cohen real. But certainly
Woodin's argument is the inspiration for both my fragile measurability
result \cite[H] and also this paper.

Because in the inner models like $L[\mu]$ the large cardinal
property is fragile and superdestructible, all these theorems---Woodin's theorem, my Fragile Measurability theorem,
and the Superdestruction theorem---tend to show how one may obtain
inner-model-like properties by forcing. For superdestructibility
this is interesting; it has the consequence that
large cardinals, in principle, cannot automatically have any amount of indestructibility.

Before beginning the proof, I would like to point out that in
response to Apter's question
Saharon Shelah has proved, independently, that small forcing makes $\k$
destructible. His technique is to code the small generic $g$ into the
continuum function above $\k$. If $\l$ is above all this coding,
then a reflection argument shows that since the
continuum function below $\k$ cannot code the new set, $\k$ cannot be
still $\l$-supercompact. Since it relies, however, on building a
particular $\ltk$-closed poset which will destroy the supercompactness
of $\k$, this technique seems not to show superdestructibility. My
argument establishes the stronger result that essentially {\it all} such
posets kill the supercompactness of $\k$.

Let's now begin my proof. I will rely on the following fact. Woodin
based the theorem I mentioned above on a similar fact concerning his
reverse Easton $\k$-iterations.

\lemma Key Lemma. If $\card{\P}=\b$, $\force_{\P}\hbox{$\Qdot$ is
$\lteb$-closed}$, and $\cof(\l)>\b$, then $\P*\Qdot$ adds no
new subset of $\l$ all of whose initial segments are in the ground
model $V$.

\proof Such sets, which are not in $V$ but all of whose initial
segments are in $V$, I will say are {\df fresh} over $V$. If the
lemma fails for some $\P$ and $\Qdot$, then we may assume
there is a name $\t$ for the characteristic function of the fresh set,
so that
$$\PQforces\t\in 2^{\check \l}\and\t\notin\check
             V\and\forall\g<\check\l\;\t\restrict\g\in\check  V.$$
By refining to a condition if necessary, we may assume that $\P$
adds a fresh subset to some minimal $\d\leq\b$, so
$$\one\Pforces\hdot\in 2^{\check \d}\and\hdot\notin\check
V\and\forall\a<\check \d\;\hdot\restrict\a\in\check V.$$
(I will actually only use that $\Qdot$ is $\lted$-closed.)
The basic idea of this proof will be to use the small set $h$ added by $\P$ to define a path through an initial segment of the tree
of attempts to decide more and more of $\t$, using the
$\lted$-closure of $\Qdot$. Since all the initial
segments of $\t$ are in $V$ we will find a set $b$ in $V$ which reveals
to us the path determined by $h$, and this will contradict the fact
that $h$ is not in $V$.
A bit of notation: if $\pq\in\PQ$,
then let $b_\pq$ be the longest sequence $b$ such that
$\pq\forces\check  b\of\t$. Also, write $\pq\decides\t\restrict\g$
to mean that $\pq$ decides $\t\restrict\g$, i.e.,
$\pq\forces \t\restrict\g=\check b$ for some $b\in V$.
The crucial aspect of the following claim
is that the first coordinate $p$ does not vary.

\claim. There is $\pq\in\PQ$ such that whenever $\gG$ is $V$-generic
below $\pq$ then for every $\g<\l$ there is a condition $\pr\in\gG$
such that $\pr\decides\t\restrict\g$.

\proof  Let $\gG$ be $V$-generic for $\PQ$. In $\VgG$ pick for every
$\g<\l$ a condition $\pgqg\in\gG$ such
that $\pgqg\decides\t\restrict\g$. Thus $\seq<p_\g\st\g<\l>$ is a
sequence of conditions from the poset $\P$. Since $\cof(\l)>\b$, and
this is preserved by $\P$ and $\Q$,
there must be some condition $\pbar$ which is repeated cofinally. In
fact, we could have used $\pbar$ in every choice. So assume that
$\<\pbar,\qdot_\g>$ decides $\t\restrict\g$ for every $\g<\l$. This fact
must be forced by some $\pq$, where $p\lte\pbar$. Thus, any generic $g*G$
containing $\pq$ satisfies $\forall\g{<}\l\;\exists\rdot\<\pbar,\rdot>\in\gG\and
\<\pbar,\rdot>\decides\t\restrict\g$. Now replace $\pbar$ with the
stronger condition $p$ to conclude the claim.\qed

Fix $\pq$ as in the claim.

\claim. For any $\pr\lte\pq$ there are $\rdot_0$ and $\rdot_1$ such
that $\<p,\rdot_0>,\<p,\rdot_1>\lte\pr$ and
$b_{\<p,\rdot_0>}\perp b_{\<p,\rdot_1>}$.

\proof If not, then some $\pr$ fails to split in that sense. Force below
$\pr$ to obtain $V$-generic $\gG$ with $\pr\in\gG$. Because of the
splitting failure, all $b_{\ps}$ with $\ps\lte\pr$ must cohere. But
by the property of the first claim, they also decide more and more of
$\t$. Thus, $\t_{\gG}=\union\set{b_\ps\st\ps\lte\pr}$, which contradicts
that $\forces \t\notin\check  V$.\qed

Iterating the claim transfinitely, I define $\qdot_t$ by
induction on $t\in 2^{\ltd}$, so that $\qdot_{\emptyset}=\qdot$ and
\points 1. $t\of\bar t\implies \<p,\qdot_{\bar
           t}>\lte\<p,\qdot_t>$\cr
        2. $b_{\<p,\qdot_{t\concat 0}>}\perp
            b_{\<p,\qdot_{t\concat 1}>}$.\cr\par

\def\qtzero{\qdot_{t\concat 0}}
\def\qtone{\qdot_{t\concat 1}}
\def\pqt{{\<p,\qdot_t>}}
\def\pqtzero{{\<p,\qtzero>}}
\def\pqtone{{\<p,\qtone>}}
\noindent At successor stages, simply apply the claim. At limit stages,
when $\qdot_t$ is defined for all $t\ofnoteq\bar t$, then
$p\Pforces\seq<\qdot_t\st t\ofnoteq\bar t>$ is descending, and so by
the closure of $\Qdot$ we obtain $\qdot_{\bar t}$.

Now force below $\pq$ so that $\pq\in\gG$ for some $V$-generic $\gG$.
Let $h=(\hdot)_g$ be the new $\d$-sequence which was added by $\P$.
Thus every initial segment $t\ofnoteq h$ is in $V$. Let
$q_t=(\qdot_t)_g$. By condition 1 it follows that
$\seq<q_t\st t\ofnoteq h>$ is a $\d$-descending sequence in
$\Q=\Qdot_g$, and so by closure there is a condition $r$
such that $r\lte q_t$ for all $t\ofnoteq h$. Let
$b=\union_{t\ofnoteq h}b_{\pqt}$. Thus, $r\Qforces b\ofnoteq \t$,
and therefore $b\in V$. But this is impossible, since $b$ will decode
for us in $V$ the generic set $h$: by construction,
$b_{\pqt}\of b$ only when $t\of h$, since condition 2 ensures
that whenever $t\concat i$ first deviates from $h$, then
$b_{\<p,\qdot_{t\concat i}>}$ will deviate from $b$. We therefore
conclude that $h\in V$, contrary to our choice.\qed

Now I am ready to prove the main theorem in parts.

\theorem Superdestruction Theorem I. Small forcing makes any
cardinal superdestructible.

\proof It suffices to show that if $\card{\P}<\k$, and
$$\Pforces\hbox{$\Qdot$ is $\ltk$-closed, and adds a new subset of
$\k$},$$ then $\k$ is not measurable after forcing with
$\PQ$. Let's suppose this fails for some $\PQ$, and that $\VgG$ is a
forcing extension by $\PQ$ in which $\k$ is measurable.
Since $\k$ is measurable, there is an embedding
$j:V[g][G]\to N$ for some transitive $N$ with $\cp(j)=\k$.
By elementarity we may decompose $N$ into it's forcing
history and write the embedding as $j:V[g][G]\to M[g][j(G)]$ for
some transitive $M$. One should not assume that $M\of V$,
since the embedding $j$ is not neccessarily the lift of an
embedding in $V$. Nevertheless, we have the following claim:

\claim. $P(\k)^M\of V$.

\proof First note that $M_\k=V_\k$ since $\cp(j)=\k$. Now suppose
that $B\of\k$ and $B\in M$. Thus, $B\intersect\a$ is in $V$ for
every $\a<\k$, and so every initial segment of $B$ is in $V$.
It follows by the Key Lemma that $B\in V$.\qed

Let
$A\of\k$ be the new set added by $\Q$, so $A\in\VgG\setminus V[g]$.
Since $A=j(A)\intersect\k$ it follows that $A\in M[g][j(G)]$.
But the $j(G)$ forcing was ${\lt}j(\k)$-closed, and so actually
$A\in M[g]$. Therefore, $A=\Adot_g$ for some name $\Adot\in M$.
We may view $\Adot$ as a function from $\k$ to the
set of anti-chains in $\P$, and this can be coded with a
subset of $\k$. So, by the claim, $\Adot\in V$, and thus
$A=(\Adot)_g\in V[g]$. This contradicts the choice of $A$.\qed

Before going on to the improved versions of the Superdestruction
Theorem, let me just point out the following corollary.

\corollary. One can force to make every large cardinal
superdestructible.

\proof Just add a Cohen real and apply the Superdestruction Theorem.\qed

That it is
so easy to make cardinals superdestructible is surprising, since in
\cite[H] a very great effort is made to make a single supercompact
cardinal have
fragile measurability. This corollary also shows that Woodin's
entire reverse Easton iteration---the one which makes the measurability of
a supercompact cardinal $\k$ destructible by $\add(\k,1)$---can be replaced by the forcing to add a Cohen real or indeed any small forcing, with
the result that every cardinal $\k$ becomes destructible by
$\add(\k,1)$, among many other posets.

\theorem Superdestruction Theorem II. After small forcing,
any $\ltk$-closed forcing which adds a subset to $\k$ will
destroy the weak compactness of $\k$.

\proof We will follow the proof of the previous theorem, but use instead
only a weakly-compact embedding. Let $V[g][G]$, etc., be as in the
earlier proof. Now suppose only that $\k$ is weakly compact in
$V[g][G]$. Pick $\l\muchgt\k$ very large, and let $X\elesub V_\l[g][G]$
be an elementary submodel of size $\k$ with $V_\k\of X$ and
$g,G,\P,\Qdot,A\in X$. The Mostowski collapse of $X$ will be a
structure $N[g][G^*]$ of size $\k$, where $N$ is transitive.
By the weak-compactness of $\k$ there is an embedding
$j:N[g][G^*]\to M[g][j(G^*)]$ for some transitive $M$ with $\cp(j)=\k$.
Since again
by the critical point we know that $M_\k=V_\k$, it follows by the
Key Lemma that $P(\k)^M\of V$. Now argue again that $A\in M[g]$ and
so $A=\Adot_g$ for some name $\Adot\in M$. But again $\Adot$ can
be thought of as a function from $\k$ to the antichains of $\P$, and
so it may be coded as a subset of $\k$. Thus, again $\Adot\in V$, and
so $A=\Adot_g\in V[g]$, contrary to the choice of $A$.\qed

Next, I push the previous arguments up to the case where the new sets
are added by $\Q$ perhaps only above $\k$.

\theorem Superdestruction Theorem III. After small forcing any cardinal
$\k$ becomes superdestructible at $\k$, at $\k\plus$, at $\k\plusplus$,
etc. In fact, if the small forcing is $\ltd$-distributive, then $\k$
becomes superdestructible at every $\l$ below $\aleph_{\k+\d}$.

\proof Suppose that $\card{\P}<\k$, that $\P$ is $\ltd$-distributive,
but, using the notation of the previous proofs, that $\k$ remains $\l$-supercompact in $V[g][G]$, where $G\of\Q$ adds a new subset
$A\of\l$, and $\l=\aleph_{\k+\b}$ for some $\b<\d$. We may assume that
$\Q$ adds no new subsets
of any smaller ordinal. In $V[g][G]$ there is a $\l$-supercompact
embedding $j:V[g][G]\to M[g][j(G)]$.

\claim. $P(\l)^M\of V$.

\proof I will show by induction that $P(\aleph_{\k+\a})^M\of V$ for all
$\a\leq\b$. To begin, we know by the argument in
the previous theorems that $P(\k)^M\of V$, since there are no new
subsets of $\k$ in $V[g][G]$ all of whose initial segments are in $V$.
Also, I claim that $(\aleph_{\k+\a})^M =(\aleph_{\k+\a})^{M[g][j(G)]}=(\aleph_{\k+\a})^{V[g][G]}=
(\aleph_{\k+\a})^V$. The first
equality holds because of the smallness of $\P$ and the
closure of $j(\Q)$. The
second equality holds because of the closure of the embedding $j$. The
last equality holds by the smallness of $\P$ and the minimality of $\l$.
Now suppose that
$P(\aleph_{\k+\a})^M\of V$, that $B\of \aleph_{\k+(\a+1)}$, and that
$B\in M$.
Every initial segment of $B$ is coded with a subset of $\aleph_{\k+\a}$
in $M$, and therefore lies in $V$ by the induction hypothesis. Since
$\aleph_{\k+(\a+1)}$ is regular, it follows by the Key Lemma that $B$ is in $V$. This completes the successor stage. Now suppose that
$P(\aleph_{\k+\a})^M\of V$ for all $\a<\g$ where $\g\leq\b$ is a limit
ordinal. If $B\of\aleph_{\k+\g}$ and $B\in M$ then again every initial segment
of $B$ is in $V$ by the induction hypothesis. But the forcing $\P$ is
$\ltd$-distributive, and $\g\leq\b<\d$, so $\P$ cannot add $B$ (this
is where the limitation on $\l$ is used).
Similarly, the highly closed $\Q$ cannot add $B$, so it must be that
$B\in V$. This establishes the limit case, and so the claim is proved.\qed

Since the embedding is closed
under $\l$-sequences, it follows that $A\in M[g][j(G)]$. But $j(G)$ is
${\lt}j(\k)$-closed and $\l<j(\k)$, so $A\in M[g]$, and thus $A=\Adot_g$ for some
name $\Adot\in M$. Again, we may view $\Adot$ as a function from
$\l$ to the set of antichains of $\P$. Thus, $\Adot$ may be coded with
a subset of $\l$ in $M$. By the claim it follows that $\Adot\in V$,
and so $A=\Adot_g\in V[g]$, contrary to our choice of $A$.\qed

\theorem Superdestruction Theorem IV. Suppose that
$\P$ has cardinality $\b$, adds a new subset to $\bar\d$, and is $\ltd$-distributive. Suppose also that $\b<\k\leq\l<\aleph_{\k+\d}$. Then any further ${\lte}\bar\d$-closed forcing which preserves
$2^{\ltk}$ but adds a subset to $\l$ will destroy the $\l$-supercompactness
of $\k$.

\proof Just apply the full power of the Key Lemma to the previous
proofs. We never used full $\ltk$-closure---rather, we used
${\lte}\bar\d$-closure to apply the Key Lemma, and we used the
preservation of $j(2^\ltk)$ by $j(\Q)$ to know that $A\in M[g]$. So the
proofs go through for the broader class of posets in this theorem.\qed

This last version of the Superdestruction Theorem is actually an
enormous improvement, reducing $\ltk$-closure to something much
less. If, for example, $\P$ is the forcing to add a Cohen real, then we
obtain the following corollary.

\corollary. After adding a Cohen real, the measurability of
any cardinal $\k$ is destroyed by any countably-closed poset which
adds a new subset, but no bounded subset, to $\k$. Similarly, the
$\l$-supercompactness of $\k$ is destroyed by any countably-closed
poset which adds a new subset to $\l$, but no bounded subset to $\k$, for
$\l=\k,\k\plus,\k\plusplus$, etc.

Finally, I will show that indestructible cardinals are not too severely
wounded when they are made superdestructible; they remain
resurrectible (this was proved, independently, by James Cummings).
My proof uses the instrumental
Term Forcing Lemma, a part of mathematical folklore,
which allows us in a sense to reverse the order of an iteration $\PQ$.

\lemma Term Forcing Lemma. If $\P*\Qdot$ is a forcing iteration, then
there is a poset $\Qterm$ such that forcing with the product $\Qterm\cross\P$ produces canonically a generic for the poset $\PQ$.
Hence, forcing with $\Qterm\cross\P$ is equivalent to forcing with
$\PQ*\Rdd$ for some (name of a) poset $\Rdd$. Finally,
if $\one\Pforces\hbox{$\Qdot$ is $\ltk$-directed closed}$, then
$\Qterm$ is also $\ltk$-directed closed.

\proof We may assume, by using a better name if necessary, that $\Qdot$
is a full name, in the sense that  if $\one\Pforces\s\in\Qdot$ then
there is a name $\t\in\dom(\Qdot)$ such that $\one\Pforces\s=\t$.
Now let $\Qterm=\set{\s\in\dom(\Qdot)\st\one\forces\s\in\Qdot}$. Define
the order $\s\termlte\t$ iff $\one\Pforces\s\Qlte\t$. Now suppose that
$\Gterm\of\Qterm$ is $V$-generic, and $g\of\P$ is $V[\Gterm]$-generic.
We must find in $V[\Gterm][g]$ a generic for $\PQ$. Let
$G=\set{\s_g\st\s\in\Gterm}$.

\claim. $g*G$ is $V$-generic for $\PQ$.

\proof We know that $g\of\P$ is $V$-generic, so it suffices to show
that $G$ is $V[g]$-generic for $\Q=\Qdot_g$. First observe that
$G$ is truly a filter, since if $\s_g$, $\t_g$ are in $G$, with
$\s,\t\in\Gterm$, then there must be some term $\eta\in\Gterm$ such
that $\eta\termlte\s,\t$. It follows that $\eta_g\Qlte\s_g,\t_g$. So
$G$ is a filter. Let's now check the genericity criterion. Suppose
that $D\of\Q$ is dense, where $D=\Ddot_g$ for some name $\Ddot$. We
may assume that $\one\Pforces\Ddot\hbox{ is dense in }\Qdot$. Now
let $\Dterm=\set{\s\in\Qterm\st\one\Pforces\s\in\Ddot}$. Observe
that $\Dterm$ is dense in $\Qterm$ since given any name $\s\in\Qterm$
we may find a name $\t$ such that
$\one\Pforces\t\Qlte\s\and\t\in\Ddot$. Thus there is a name
$\s\in\Gterm\intersect\Dterm$, and so $\s_g\in G\intersect D$.\qed

Since $\Qterm\cross\P$ produces a generic $g*G$ for $\PQ$, it follows that
the regular open algebra of $\PQ$ completely embeds into the regular
open algebra of $\Qterm\cross\P$ via the map
$$\pq\mapsto\boolval{\pq\in g*G}^{\Qterm\cross\P}.$$ By
standard quotient forcing arguments (see, e.g.,\cite[J] p. 237, ex.
23.6), it follows that forcing with $\Qterm\cross\P$
is equivalent to forcing with $\PQ*\Rdd$ for some (name of a) poset
$\Rdd$.

It remains to prove the last sentence of the lemma. Suppose that
$$\one\Pforces\Qdot\hbox{ is $\ltk$-directed closed},$$ and that
$A\of\Qterm$ is a $\ltk$ size family with the FIP. With a slight
abuse of name notation, it follows that $\one\Pforces\hbox{$A\of\Qdot$
is a $\ltk$ size family with the FIP}$. Using the directed closure of
$\Qdot$, we obtain a name $\s$ such that $\one\forces\s\Qlte\t$ for
every $\t\in A$. Thus, $\s\termlte\t$ for every $\t\in A$, and the
lemma is proved.\qed

\theorem Resurrection Theorem. After small forcing an indestructible
cardinal remains resurrectible, but never strongly resurrectible.

\proof The implicit claim of this theorem, that
indestructible cardinals are resurrectible, is clear: if $\k$
is indestructible, and $\Q$ is $\ltk$-directed closed, then
$\k$ is supercompact in $V^{\Q}$. So no further forcing needs to
be done to recover the supercompactness of $\k$. Thus, indestructible
cardinals are in fact strongly resurrectible.

Now suppose that $\P$ is small. I will show that $\k$ remains
resurrectible after forcing with $\P$. So suppose
$\Pforces\Qdot\hbox{ is $\ltk$-directed closed}$. I want to
recover the supercompactness of $\k$ by further forcing after
$\PQ$. By the
Term Forcing Lemma, forcing with
$\Qterm\cross\P$ is equivalent to forcing with $\PQ*\Rdd$, for
some $\Rdd$.
Furthermore, $\Qterm$ is $\ltk$-directed closed
and therefore preserves the supercompactness of $\k$,
since $\k$ was indestructible in $V$. Small forcing by $\P$ then
also preserves the supercompactness of $\k$. Thus, forcing with
$\Qterm*\P$, and hence also $\PQ*\Rdd$, preserves the supercompactness of
$\k$. Therefore, the forcing $\Rdd$ over $V^{\PQ}$ must have recovered the
supercompactness of $\k$.

It remains to check that $\Rdd$ is sufficiently distributive. That is,
we have to show that $\Rdd$ adds no new $\g$-sequences for any $\g<\k$.
Suppose
that $\Gterm*g\of\Qterm*\P$ is $V$-generic, and produced the generics
$g*G*H\of\PQ*\Rdd$, where $V[\Gterm][g]=V[g][G][H]$. Suppose that
$s\in V[\Gterm][g]$ is a $\g$-sequence of ordinals for some $\g<\k$.
So $s=\sdot_g$ for some $\sdot\in V[\Gterm]$, where $\sdot$ is a
function from $\g$ to a (small) set of antichains in $\P$ matched with
ordinals (i.e. the possible values of $\sdot(\g)$). It follows that
$\sdot\in V$ since $\Qterm$ is $\ltk$-directed closed. Thus, $s\in V[g]$. Thus, the only $\g$-sequences added by $H$ must have been
already in $V[g]$, and so $\Rdd_{g*G}$ must be $\ltk$-distributive.
So $\k$ is resurrectible in $V[g]$.

Finally, I will
show that $\k$ is not strongly resurrectible in $V[g]$.
Let $\Q$ be the poset in $V[g]$ to add a Cohen subset to
$\k$, or in fact any $\ltk$-closed poset which adds a subset to $\k$.
We know by the Superdestruction Theorem that $\Q$ will destroy
the measurability of $\k$. If $\Rdot$ is the $\Q$-name of a
$\ltk$-closed poset in $V[g]^\Q$, then it follows that $\Q*\Rdot$ is
$\ltk$-closed in $V[g]$, since the $*$-iteration of closed posets is
closed. Since it also
adds a subset to $\k$ it follows again by the Superdestruction Theorem
that $\Q*\Rdot$ will destroy the measurability of $\k$. Thus, the
supercompactness of $\k$ cannot be recovered by $\ltk$-closed
forcing. Therefore, $\k$ is not strongly resurrectible in $V[g]$.\qed


Let me list, finally, two natural questions which remain unanswered in this
paper. The first asks whether the limitation on $\l$ in the Superdestruction Theorem III can be removed. The second asks more generally
whether small forcing leads to a certain attractive complement of Laver
indestructibility.

\question. After small forcing, does $\k$ become superdestructible at
$\l$ for every $\l$?

\question. After small forcing, does every $\ltk$-closed forcing destroy
the supercompactness of $\k$?

Though I have not answered these questions in this paper, I nevertheless
know that the answer to both of them is `yes'. In a forthcoming
paper which I am now writing with Saharon Shelah, we prove that
after small forcing, any $\ltk$-closed
forcing will destroy even the strong compactness of $\k$. Thus, after
small forcing, a supercompact cardinal $\k$ has a dual version of
Laver indestructibility. Namely, it is destroyed by any $\ltk$-closed
forcing.

\bigskip
{\tenpoint\noindent
I would like to thank AnnMarie Fela at Fela's Cafe, now
health-consciously reincarnated as The Secret Garden,
for making me such delicious pancakes while I proved the theorems in this
paper.}

\section Bibliography

\tenpoint
\nopagenumbers
\parindent=0pt
\newbox\Article
\newbox\Journal
\newbox\Author
\newbox\Vol
\newbox\No
\newbox\Year
\newbox\Page
\newbox\Book
\newbox\Publisher
\newbox\Pubaddr
\newbox\Key
\newbox\Editor
\newbox\Comment
\newbox\Note
\def\entry#1#2\par{\item{#1\quad}\hskip-1.1em#2\par}
\def\article#1{\setbox\Article=\hbox{\sl #1, }}
\def\journal#1{\setbox\Journal=\hbox{\rm #1 }}
\def\author#1{\setbox\Author=\hbox{\sc #1, }}
\def\vol#1{\setbox\Vol=\hbox{\bf #1 }}
\def\no#1{\setbox\No=\hbox{no. #1 }}
\def\year#1{\setbox\Year=\hbox{\rm({\oldstyle #1}) }}
\def\page#1{\setbox\Page=\hbox{\rm p. #1 }}
\def\book#1{\setbox\Book=\hbox{\it #1, }}
\def\publisher#1{\setbox\Publisher=\hbox{\rm #1, }}
\def\pubaddr#1{\setbox\Pubaddr=\hbox{\rm #1, }}
\def\key#1{\setbox\Key=\hbox{#1}}
\def\editor#1{\setbox\Editor=\hbox{\rm(#1, Ed.) }}
\def\comment#1{\setbox\Comment=\hbox{\rm #1}}
\def\note#1{\setbox\Note=\hbox{\rm #1 }}
\def\ref#1\par{\smallskip{#1
  \entry{\ifhbox\Key\unhbox\Key\else[\ ]\fi}%
  \unhbox\Author\unhbox\Note
  \ifhbox\Book \unhbox\Book\unhbox\Publisher\unhbox\Pubaddr
               \unhbox\Editor\unhbox\Page\unhbox\Year\unhbox\Comment
  \else \unhbox\Article\unhbox\Journal\unhbox\Vol\unhbox\No\unhbox\Editor
        \unhbox\Page\unhbox\Year\unhbox\Comment\fi\par}}
  
\ref
\article{Menas' Result Is Best Possible}
\author{Arthur Apter \& Saharon Shelah}
\journal{to appear in Transactions of the AMS}
\key{[AS]}

\ref
\article{Making the hugeness of $\k$ resurrectable after 
$\k$-directed closed forcing}
\author{Julius B. Barbenel}
\journal{Fundamenta Mathematicae}
\vol{137}
\year{1991}
\page{9-24}
\key{[B]}

\ref
\author{Joel David Hamkins}
\article{Fragile Measurability}
\journal{Journal of Symbolic Logic}
\year{1994}
\vol{59}
\page{262-282}
\key{[H]}

\ref
\author{Thomas Jech}
\book{Set Theory}
\publisher{Academic Press}
\year{1978}
\pubaddr{London}
\key{[J]}

\ref
\author{Richard Laver}
\article{ Making the Supercompactness of $\kappa$ Indestructible Under 
 $\kappa$-Directed Closed Forcing}
\journal{Isreal Journal Math}
\vol{29}
\year{1978}
\page{385-388}
\key{[L]}

\ref
\author{W. Hugh Woodin}
\article{Forcing to a Model of a Supercompact $\k$ Whose Weak Compactness
is Killed by $\add(\k,1)$}
\comment{unpublished theorem}
\key{[W]}

\bye